\documentclass[12pt, reqno]{amsart}

\usepackage{amsthm, amssymb}
\usepackage{txfonts,exscale}

\newtheorem{thm}{Theorem}
\newtheorem*{lem}{Lemma}

\newtheorem{question}{Question}

\theoremstyle{definition}

\newtheorem*{acknowledgement}{Acknowledgement}

\title{Approximating reals by sums of rationals}
\author[T. H. Chan]{Tsz Ho Chan}
\address{American Institute of Mathematics\\360 Portage Avenue\\Palo Alto, CA 94306\\U.S.A.}
\email{thchan@aimath.org}
\author[A.V. Kumchev]{Angel V. Kumchev}
\address{Department of Mathematics \\1 University Station, C1200 \\The University of Texas at Austin \\Austin, TX 78712 \\U.S.A.}
\email{kumchev@math.utexas.edu}

\begin{document}

\begin{abstract}
  We study how well a real number can be approximated by sums of two or more rational numbers with denominators up to a certain size.
\end{abstract}

\maketitle

\section{Introduction and main result}

Dirichlet's theorem on diophantine approximation tells us that we can approximate any real number by rational numbers quite well, namely:

\begin{thm}\label{theorem1}
  For any real $\theta$ and any positive integer $N$, there exist integers $a$ and $q$, with $1 \leq q \leq N$, such that
  \[
    \left| \theta - \frac{a}{q} \right| < \frac {1}{qN}.
  \]
\end{thm}

Moreover, the bound $1/(qN)$ is best possible, apart from the constant factor. To see this, it suffices to consider the golden ratio $\theta = (\sqrt{5} - 1)/2$ (see \cite[\S 11.8]{HW}). During his work in \cite{C}, the first author accidentally stumbled across the following analogous question:

\begin{question}\label{q1}
  For any real $\theta$ and any positive integer $N$, give an upper bound for 
  \[
    \min_{ \substack{ a_1, a_2, q_1, q_2 \in \mathbb{Z}\\ 1 \leq q_1, q_2 \leq N}} \left| \frac{a_1}{q_1} + \frac{a_2}{q_2} - \theta \right|.
  \]
\end{question}

With the golden ratio in mind, we know that the upper bound can be no better than $O\big( 1/(q_1 q_2 N^2) \big)$. So, what is the best possible upper bound? More generally,

\begin{question}\label{q2}
  Let $k$ be a positive integer. For any real $\theta$ and any positive integer $N$, give an upper bound for
  \[
    \min_{ \substack{ a_1, \dots, a_k, q_1, \dots, q_k \in \mathbb{Z}\\ 1 \leq q_1, \dots, q_k \leq N}} \left| \frac{a_1}{q_1} + \dots + \frac{a_k}{q_k} - \theta \right|.
  \]
\end{question}

To these, we have the following result:

\begin{thm}\label{theorem2}
  Let $k$ be a positive integer. For any real $\theta$ and any positive integer $N$, there exist integers $a_1, \dots, a_k$, $q_1, \dots, q_k$, with $1 \le q_1, \dots q_k \le N$, such that
  \[ 
    \left| \frac {a_1}{q_1} + \dots + \frac {a_k}{q_k} - \theta \right| \ll N^{-k}. 
  \]
\end{thm}

  The bound $N^{-k}$ is best possible in the sense that, for some $\theta$, the minimum in Question \ref{q2} can be as large as $N^{-k}$. For example, if one considers $\theta = 1/(2N^k)$,
  \[ 
    \left| \frac{a_1}{q_1} + \cdots + \frac{a_k}{q_k} - \theta \right| \ge \frac 1{2N^k}
  \]
  for any choice of $a_1, \dots, a_k, q_1, \dots, q_k$, with $1 \le q_1, \dots, q_k \le N$. However, one expects such pathological examples to be relatively rare, and so one may wonder if it is possible to obtain a sharper upper bound involving the $q_i$'s. For example, is it possible to replace $N^{-k}$ by $(q_1 \cdots q_k)^{-1}N^{-k}$ in Theorem \ref{theorem2}? We shall briefly address this issue in the last section.

\section{Proof of Theorem \ref{theorem2}}
\begin{lem}
\label{lemma1}
Suppose that $k \ge 1$ is an integer. There is a number $x_0(k) \ge 1$ such that
\[ \sideset{}{^*} \sum_{1 \le q_1, \dots, q_k \le x} q_1 \cdots q_k \gg x^{2k}, \] whenever $x \ge x_0(k)$. Here, $\sum^*$ denotes a summation over the $k$-tuples $q_1, \dots, q_k$ such that $(q_i, q_j) = 1$ whenever $1 \le i < j \le k$.
\end{lem}

\begin{proof}
  It suffices to show that
  \begin{equation}\label{1}
    \sum_{ \substack{ n \le x\\ (n, m) = 1}} \phi(n)^{\alpha}n^{1 - \alpha} \gg x^2\phi(m)m^{-1},
  \end{equation}
  whenever $0 \le \alpha \le k - 1$, $1 \le m \le x^{k - 1}$, and $x \ge x_0(k)$. The conclusion of the lemma will then follow by successive applications of \eqref{1} with $\alpha = 0, 1, \dots, k - 1$ to the summations over $q_k, q_{k - 1}, \dots, q_1$.

  We now proceed to establish \eqref{1}. We start by showing that
  \begin{equation}\label{2}
    \sum_{ \substack{ n \le x\\ (n, m) = 1}} \left( \frac n{\phi(n)} \right)^{\alpha} \ll x \phi(m)m^{-1}.
  \end{equation}
  Define the multiplicative functions
  \[ 
    f(n) = \begin{cases}
      \big( n/\phi(n) \big)^{\alpha} & \text{if } (n, m) = 1, \\
      0 & \text{if } (n, m) > 1,
    \end{cases} \qquad g(n) = \sum_{d \mid n} f(d)\mu(n/d). 
  \]
  Then $g(n) \ge 0$, and
  \begin{align*}
    \sum_{n \le x} f(n) &= \sum_{n \le x} \sum_{d \mid n} g(d) = \sum_{d \le x} g(d) \left\lfloor \frac xd \right\rfloor \le x \sum_{d \le x} g(d)d^{-1}\\
    &\le x \prod_{p \le x} \sum_{\nu = 0}^{\infty} g(p^{\nu})p^{-\nu} = x \prod_{p \le x} \left( 1 - p^{-1} \right) \sum_{\nu = 0}^{\infty} f(p^{\nu})p^{-\nu} \\
    &\le x \prod_{ \substack{ p \mid m\\ p \le x}} \left( 1 - p^{-1} \right) \prod_p \left( 1 + \frac {p^{\alpha} - (p - 1)^{\alpha}}{p(p - 1)^{\alpha}} \right) \\
    &\le x \prod_{p \mid m} \left( 1 - p^{-1} \right) \prod_p \left( 1 + \frac {p^{\alpha} - (p - 1)^{\alpha}}{p(p - 1)^{\alpha}} \right) + O(1). 
  \end{align*}
  The last inequality follows on noting that $m$ has at most $k - 2$ prime divisors $p > x$, and hence,
  \[
    \prod_{ \substack{ p \mid m\\ p > x}} \left( 1 - p^{-1} \right) = 1 + O \big( x^{-1} \big). 
  \]
  This proves \eqref{2}. On the other hand, when $\alpha = 0$, we have
  \[ 
    \sum_{ \substack{ n \le x\\ (n, m) = 1}} n = \sum_{d \mid m} \mu(d) d \sum_{k \le x/d} k = \frac {\phi(m)}{2m} x^2 + O(x\tau(m)), 
  \]
  whence
  \begin{equation}\label{3}
    \sum_{ \substack{ n \le x\\ (n, m) = 1}} n^{1/2} \ge x^{-1/2} \sum_{ \substack{ n \le x\\ (n, m) = 1}} n \gg x^{3/2}\phi(m)m^{-1}.
  \end{equation}
  Finally, \eqref{1} follows from \eqref{2}, \eqref{3}, and Cauchy's inequality:
  \[ 
    \sum_{ \substack{ n \le x\\ (n, m) = 1}} \phi(n)^{\alpha}n^{1 - \alpha} \ge \bigg\{ \sum_{ \substack{ n \le x\\ (n, m) = 1}} n^{1/2} \bigg\}^2 \bigg\{ \sum_{ \substack{ n \le x\\ (n, m) = 1}} \left( \frac n{\phi(n)} \right)^{\alpha} \bigg\}^{-1} \gg x^2\phi(m)m^{-1}. 
  \]
\end{proof}

\begin{proof}[Proof of Theorem \ref{theorem2}]
  For $0 < \Delta < 1/2$, define
  \[ 
    t(x) = \max \big( 1 - |x|/\Delta, 0 \big), \qquad g(x) = \sum_{n = -\infty}^{\infty} t(x - n). 
  \]
  The function $g$ has a Fourier expansion
  \[ 
    g(x) = \sum_{h = -\infty}^{\infty} \hat g_h e(h x), \qquad \hat g_h = \Delta \left( \frac{\sin \pi \Delta h}{\pi \Delta h} \right)^2. 
  \]
  We consider the sum
  \begin{equation}\label{4}
    \mathcal S = \sideset{}{^*} \sum_{1 \le q_1, \dots, q_k \le N} \sum_{a_1 = 1}^{q_1} \dots \sum_{a_k = 1}^{q_k} g \left( \frac{a_1}{q_1} + \dots + \frac{a_k}{q_k} - \theta \right),
  \end{equation}
  where $\sum^*$ has the same meaning as in the Lemma. Putting in the Fourier expansion for $g$, we get
  \begin{align}\label{5}
    \mathcal S =& \sideset{}{^*} \sum_{1 \le q_1, \dots, q_k \le N} \sum_{a_1 = 1}^{q_1} \dots \sum_{a_k = 1}^{q_k} \sum_{h = -\infty}^{\infty} \hat g_h e \left( h \left ( \frac{a_1}{q_1} + \dots + \frac {a_k}{q_k} - \theta \right) \right) \\
    =& \sideset{}{^*} \sum_{1 \le q_1, \dots, q_k \le N} \sum_{h = -\infty}^{\infty} \hat g_h e(-h \theta) \sum_{a_1 = 1}^{q_1} e \big( ha_1/ q_1 \big) \dots \sum_{a_k = 1}^{q_k} e \big( ha_k/ q_k \big) \notag\\
    =& \sideset{}{^*} \sum_{1 \le q_1, \dots, q_k \le N} q_1 \cdots q_k \sum_{ \substack{ h = -\infty\\ q_1 \cdots q_k \mid h}}^{\infty} \hat g_h e(-h \theta) \notag,
  \end{align}
  as
  \[ 
    \sum_{a = 1}^{q} e(ha/q) = \begin{cases}
      q & \text{if } q \mid h, \\
      0 & \text{otherwise}.
    \end{cases} 
  \]
  If $m$ is a positive integer and $\Delta \le m^{-1}$, we have 
  \begin{equation}\label{6}
  \begin{split}
    \sum_{h \ne 0} \big| \hat g_{mh} \big| &\le 2 \bigg\{ \sum_{h = 1}^H \Delta + \frac 1{\Delta m^2} \sum_{h = H + 1}^{\infty} h^{-2} \bigg\} \\
    &\le 2 \bigg( H\Delta + \frac 1{H\Delta m^2} \bigg) \le 6m^{-1},
  \end{split}
  \end{equation}
  where $H = \big\lceil (\Delta m)^{-1} \big\rceil$; whereas if $\Delta > m^{-1}$, we have
  \begin{equation}\label{7}
    \sum_{h \ne 0} \big| \hat g_{mh} \big| \le \frac {2\zeta(2)}{\Delta m^2} \le 4m^{-1}.
  \end{equation}
  Putting \eqref{6} and \eqref{7} (with $m = q_1 \cdots q_k$) into \eqref{5}, we obtain
  \[ 
    \mathcal S = \Delta \sideset{}{^*} \sum_{1 \le q_1, \dots, q_k \le N} q_1 \cdots q_k + O \big( N^k \big),
  \]
  the $O$-implied constant being absolute (in fact, it is $6$). Therefore, upon choosing $\Delta = cN^{-k}$ with a sufficiently large $c > 0$, it follows from the Lemma that $\mathcal S > 0$. Hence, by \eqref{4},
  \[ 
    g\left( \frac {a_1}{q_1} + \dots + \frac {a_k}{q_k} - \theta \right) > 0 
  \]
  for some integers $a_1, \dots, a_k, q_1, \dots, q_k$ with $1 \le q_1, \dots, q_k \le N$. Then, by the definition of $g$,
  \[ 
    \left| \frac {a_1}{q_1} + \dots + \frac {a_k}{q_k} - n - \theta \right| \le cN^{-k} 
  \]
  for some integer $n$. This establishes the theorem.
\end{proof}

\section{Closing remarks}
We conclude this note with a short discussion of possible improvement on the bound $N^{-k}$ in Theorem \ref{theorem2}. For example, is it possible to replace $N^{-k}$ by $(q_1 \cdots q_k)^{-1} N^{-k}$? While such a result may appear to be the right generalization of Dirichlet's theorem, it is not true in general. Indeed, suppose that for any real $\theta$, there exist integers $a_1, \dots, a_k, q_1, \dots, q_k$, with $1 \le q_1, \dots, q_k \le N$, such that
\begin{equation}\label{8}
\left| \frac{a_1}{q_1} + \cdots + \frac{a_k}{q_k} - \theta \right| \le \frac C{q_1 \cdots q_kN^k}.
\end{equation}
Then
\begin{equation}\label{9}
[0, 1] \subseteq \bigcup_{q \in \mathcal D_k(N)} \bigcup_{0 \le a \le q} \left\{ \theta \in \mathbb R \; : \; |\theta - a/q| \le C/(qN^k) \right\},
\end{equation}
where $\mathcal D_k(N)$ denotes the set of least common denominators of the sums appearing on the left side of \eqref{8}. By a result of Erd\"os \cite{E}, $\mathcal D_k(N)$ has cardinality
  \[ 
    |\mathcal D_k(N)| \ll N^k(\log N)^{-c}
  \]
  for some constant $c = c(k) > 0$, so it follows from \eqref{9} that
  \[
    1 \le \sum_{q \in \mathcal D_k(N)} \sum_{0 \le a \le q} \frac {2C}{qN^k} \le 4CN^{-k}|\mathcal D_k(N)| \ll (\log N)^{-c},
  \]
  which is impossible when $N \to \infty$. On the other hand, one may hypothesize that the set of fractions with denominators in $\mathcal D_k(N)$ is distributed similarly to the set of all fractions $a/q$ with denominators $q \le N^k$. Under such a hypothesis, one might hope for an estimate with $|\mathcal D_k(N)|^{-1}$ in place of the term $(\log 3N)^cN^{-k}$ on the right side of \eqref{10} below, and such an estimate, if true, would be essentially best possible. However, upon observing that
  \[ 
    |\mathcal D_k(N)| \ge \sum_{q_1 \le N} \cdots \sum_{q_k \le N} d(q_1 \cdots q_k)^{-1} \ge \bigg\{ \sum_{q \le N} d(q)^{-1} \bigg\}^{k} \gg N^k(\log 3N)^{-k},
  \]
  we will take a more cautious approach and pose the following
\begin{question}
Let $k$ be a positive integer. Determine the least value of $c_k$ such that for any real $\theta$ and any positive integer $N$, there exist integers $a_1, \dots, a_k$, $q_1, \dots, q_k$, with $1 \le q_1, \dots, q_k \le N$, such that
\begin{equation}\label{10} 
\left| \frac {a_1}{q_1} + \dots + \frac {a_k}{q_k} - \theta \right| \ll \frac {(\log 3N)^{c_k}}{q_1 \cdots q_k N^k}.
\end{equation}
\end{question}
We leave the answer to this question to the future.

\bigskip

\begin{acknowledgement} 
  The first author would like to thank the American Institute of Mathematics for support. The second author would like to thank Jeff Vaaler for several enlightening conversations on this and related topics.
\end{acknowledgement}

\end{document}